\documentclass[a4paper,10pt]{article}
\usepackage{stmaryrd}
\usepackage{amscd}
\usepackage{bbding}
\usepackage{amsfonts}
\usepackage{cite}
\usepackage{mathrsfs}
\usepackage{graphicx,latexsym,amsmath,amssymb,amsthm,enumerate}
\usepackage{amsmath,bm}
\usepackage{color}
\theoremstyle{definition}
\newtheorem{theorem}{Theorem}[section]
\newtheorem{lemma}{Lemma}[section]

\newtheorem{problem}{Problem}[section]

\newtheorem{definition}{Definition}[section]

\hfuzz=\maxdimen
\begin{document}
\title{{\Large\bf{A new class of fractional impulsive differential hemivariational inequalities with an application}}\thanks{This work was supported by the National Natural Science Foundation of China (11471230, 11671282).}}
\author{{Yun-hua Weng$^a$, Tao Chen$^a$, Nan-jing Huang$^a$\thanks{Corresponding author.  E-mail address: nanjinghuang@hotmail.com; njhuang@scu.edu.cn}  and Donal O'Regan$^b$} \\
{\small\it a. Department of Mathematics, Sichuan University, Chengdu,
	Sichuan 610064, P.R. China}\\
{\small\it b. School of Mathematics, Statistics and Applied Mathematics, National University of Ireland, Galway, Ireland}}
\date{ }
\maketitle
\begin{flushleft}
\hrulefill\\
\end{flushleft}
{\bf Abstract.}
We consider a new fractional impulsive differential hemivariational inequality which captures the required characteristics of both the hemivariational inequality and the fractional impulsive differential equation within the same framework. By utilizing a surjectivity theorem and a fixed point theorem, we establish an existence and uniqueness theorem for such a problem.  Moreover, we investigate the perturbation problem of the fractional impulsive differential hemivariational inequality to prove a convergence result which describes the stability of the solution in relation to perturbation data. Finally, our main results are applied to obtain some new results for a frictional contact problem with the surface traction driven by the fractional impulsive differential equation.
\\ \ \\
{\bf Keywords}: Fractional differential variational inequality; Fractional impulsive equation; Hemivariational inequality; Frictional contact.
 \\ \ \\
\textbf{2020 AMS Subject Classification:} 34A08; 34G20; 49J40;  74M10; 74M15.

\section{Introduction}

Let $Y$, $Z_{1}$, $Z_{2}$ be three reflexive and separable Banach spaces and let $Z_{2} ^{*}$ be the dual space of $Z_{2}$. For a prefixed $T>0$, let $Q=[0,T]$  and
\begin{align*}
f:Q\times Z_1\times Z_2\to Z_1,\quad A:Q\times Z_2\to Z_2^{*},\quad N:Z_{2}\to Y,\quad
J:Q\times Y\to R, \quad g:Q\times Z_1\to Z_2^{*}.
\end{align*}
This paper focuses on the following fractional impulsive differential hemivariational inequality (FIDHVI): find $z:Q\to Z_1$ and $y:Q\to Z_2$ such that
$$
\begin{cases}
^{C}D_{0}^{\kappa}z(t)= f(t,z(t),y(t)),\quad   t\in Q,\;t\neq \tau_j,\;j=1,2,\cdots,m, \\
 \Lambda z(\tau_j)=\Theta_j(z(\tau_j ^{-})),\;j=1,2,\cdots,m,\\
z(0)=z_0,\\
\langle A(t,y(t)),x\rangle+J^{\circ}(t,Ny(t);Nx)\geq\langle g(t,z(t)),x\rangle,\quad \forall (t,x)\in  Q\times Z_2,
\end{cases}
$$
where $^{C}D_{0}^{\kappa}(0 < \kappa\leq 1)$ stands for the Caputo derivative of fractional order $\kappa$, $\Theta_j:Z_1\to Z_1$ is an impulsive function with $j=1,2,\cdots,m$,  $\Lambda z(\tau_j)$ is given by $\Lambda z(\tau_j)=z(\tau_j ^{+})-z(\tau_j ^{-})$ with $z(\tau_j ^{+})$ and $z(\tau_j ^{-})$ being the left and right limit of $z$ at $t=\tau_j$, respectively, and $0=\tau_0<\tau_1<\cdots<\tau_m<\tau_{m+1}=T$.

Some particular cases of FIDHVI are as follows.
\begin{itemize}
\item [(i)] If $\kappa=1$, then FIDHVI turns into  impulsive differential hemivariational inequality:  find $z:Q\to Z_1$ and $y:Q\to Z_2$ such that
$$
\begin{cases}
z'(t)= f(t,z(t),y(t)),\quad   t\in Q,\;t\neq \tau_j,\;j=1,2,\cdots,m, \\
\Lambda z(\tau_j)=\Theta_j(z(\tau_j ^{-})),\;j=1,2,\cdots,m,\\
z(0)=z_0,\\
\langle A(t,y(t)),x\rangle+J^{\circ}(t,Ny(t);Nx)\geq\langle g(t,z(t)),x\rangle,\quad \forall (t,x)\in  Q\times Z_2,
\end{cases}
$$
which is still a new problem.

\item [(ii)] If we set $\Theta_j=0$, $J=0$ and replace $f(t,z(t),y(t))$, $A(t,y(t))$ and $g(t,z(t))$ by $Bz(t)+Cy(t)$, $Az(t)$, and $-g(t,z(t))$, respectively, where $B:Z_1\to Z_1^{*}$ , $C:Z_2\to Z_1$ are given mappings,  moreover, the solutions of the hemivariational inequality are governed by  a constraint set $K\subset Z_2$,  then FIDHVI reduces to the following DVI
$$
\begin{cases}
z'(t)= Bz(t)+Cy(t),\quad   t\in Q, \\
z(0)=z_0,\\
\langle Ay(t)+g(t,z(t)),x-y(t)\rangle\geq0,\quad \forall (t,x)\in  Q\times K(\subset Z_2),
\end{cases}
$$
which was discussed by Guo et al. in \cite{zhtj9}.

\item [(iii)] If we set $\Theta_j=0$, $N=I$ and replace $f(t,z(t),y(t))$, $A(t,y(t))$, $g(t,z(t))$ and $J(t,y(t))$ by $Bz(t)+F(z(t),y(t))$, $Az(t)$,  $g(z(t),y(t))$ and $Jy(t)$, respectively, where $B:Z_1\to Z_1^{*}$ , $F:Z_1\to Z_1^{*}$ are given mappings and $I:Z_2\to Z_2$ is the identity operator, moreover, the solutions of the hemivariational inequality are governed by  a constraint set $K\subset Z_2$,  then FIDHVI takes the following form
$$
\begin{cases}
^{C}D_{0}^{\kappa}z(t)= Bz(t)+F(z(t),y(t)),\quad   t\in Q, \\
z(0)=z_0,\\
\langle Ay(t)-g(z(t),y(t)),x-y(t)\rangle+J^{\circ}(y(t);x-y(t))\geq0,\quad \forall (t,x)\in  Q\times K(\subset Z_2),
\end{cases}
$$
which was considered by Jiang et al. in \cite{attractor}.
\end{itemize}

We remark that for appropriate and suitable choices of the spaces and the above defined maps, FIDHVI
includes a number of fractional differential hemivariational inequalities, impulsive differential hemivariational inequalities, fractional differential variational inequalities, impulsive differential variational inequalities, differential hemivariational inequalities, and differential variational inequalities as special cases, see for example \cite{LXS2016,Xuegm,special1,Zeng2,MS,Liu4} and the related references cited therein.

It is worth  mentioning that FIDHVI is a new model which captures the  required characteristics of both the hemivariational inequality and the fractional impulsive differential equation within the same framework.  In addition,  FIDHVI can be used to describe the frictional contact problem with the surface traction driven by the fractional impulsive differential equation (see Section 5).

The study of  differential variational inequality (DVI)  can  ascend  to the work of Aubin and Cellina \cite{Aubin}.  DVI described by the following generalized abstract system
$$
\begin{cases}
\dot{y}(t)=f(t,y(t),x(t)),\quad \forall t\in Q,\\
 \int_{0} ^{t}\left(v-x(t)\right)^{T}F(t,y(t),x(t))dt\geq 0, \quad \forall v\in K, \\
G(y(0),y(T))=0
\end{cases}
$$
was then examined by Pang and Stewart \cite{Pang} in finite dimension Euclidean spaces.  Here $K$ is a nonempty, closed, and convex subset of $R^m$, $f:Q\times R^{n}\times R^{m}\to R^{n}$, $F:Q\times R^{n}\times R^{m}\to R^{m}$, and $G:R^{n}\times R^{n}\to R^{n}$ are three given functions.  As pointed out by Pang and Stewart \cite{Pang}, DVI provides a powerful tool of describing many practical problems such as fluid mechanical problems, engineering operation research, dynamic traffic networks, economical dynamics and frictional contact problems (\cite{C1,C2,C4,C5,C6,C7}). In 2010, Li et al. \cite{LiNJ}  discussed the solvability for a class of DVI in finite dimensional spaces.  Later, Chen and Wang \cite{CXWZ} employed the regularized time-stepping method to consider a class of parametric DVI, and provided convergence analysis for this method in finite dimensional Euclidean spaces. Liu et al.\cite{Liu2} studied  a class of nonlocal semilinear evolution DVI in Banach spaces. By using the  theory of topological degree, they  obtained some existence results for their model under some suitable assumptions.  Recently, in order to describe a free boundary problem raising from  contact mechanics, Sofonea et al. \cite{Sof2020} studied  a differential quasivariational inequality and proved the stability of the solutions for such a problem. For more works related to DVIs, we refer the reader to \cite{zhtj1,zhtj2,zhtj3,zhtj4,zhtj5,zhtj7,zhtj8,szllyy,Liu1,Liu3,Liu4,M+1,Wu} and the the references therein.

  As is well known, fractional calculus, that is, the noninteger calculus, allows us to define derivatives of arbitrary order and  has many applications  in
  practical problems \cite{FDE1,c17,c18}. Recently, by applying the fixed point approach,  Ke et al. \cite{KeTD2015} discussed the solvability of a  class of  fractional DVI  in finite dimensional spaces. Using the Rothe method, Zeng et al. \cite{Zeng1} studied a class of parabolic fractional differential hemivariational inequalities in Banach spaces. Mig\'{o}rski and Zeng \cite{MS} considered another class of fractional DVI in Banach spaces by the discrete approximation method.  Xue et al. \cite{Xuegm} discussed the existence of the mild solutions of a class of fractional DVIs in Banach spaces under some appropriate hypotheses. Very recently, Weng et al. \cite{Weng} considered a fractional nonlinear evolutionary delay system driven by a hemi-variational inequality in Banach spaces and established an existence theorem for such a system by employing the KKM theorem, fixed point theorem for condensing set-valued operators and the theory of fractional calculus. Weng et al. \cite{Weng+}
 introduced a new fractional nonlinear system described by a fractional nonlinear differential equation and a quasi-hemivariational inequality in Banach spaces and obtained some results concerned with the existence, uniqueness and the stability of the solution for such a problem under mild assumptions.

It is worth noting that, in the real world, many systems are often disturbed suddenly, and systems changes suddenly in a short time. These phenomena are called impulsive effects. We note that  diverse numerical methods and theoretical results have been widely studied for differential equations with impulsive effects  using different assumptions in the literature; for instance, we refer the reader to \cite{mc10,mcwffcbc2,mc20,mc30,mc40,mcwffcbc1} and the references therein. Recently, Li et al. \cite{LXS2016} introduced a class of impulsive DVI in finite dimensional spaces and presented some existence and stability results of the solutions under some suitable assumptions.  However, in some practical situations applications, it is necessary to consider FIDHVI. To illustrate this  point, a fractional contact problem with the surface traction driven by the fractional impulsive differential equation will be considered as an application of FIDHVI in Section 5.  The discipline of FIDHVI is still not explored and very little is known. To fill this gap, in this paper, we seek to make a contribution in this new direction.

The outline  of this work is  as follows. In the next section we  present some necessary preliminaries and notations.  After that Section 3 establishes an  existence and uniqueness result concerning FIDHVI under some mild conditions. In Section 4 we provide a stability result of the  solution of FIDHVI with  respect to the perturbation of data. Finally, we apply our main results for FIDHVI to the frictional contact problem with the surface traction driven by the fractional impulsive differential equation in Section 5.

\section{Preliminaries}
\setcounter{equation}{0}
For a Banach space $X$,  we use the following notations:
\begin{flalign*}
&C(Q;X): \; \mbox{The space of all  functions $x:Q\to X$ that is continuous};\\
&L^{p}(Q;X): \; \mbox{The space of all $p$-th power  Bochner integrable functions on $Q$ taking values in $X$};\\
&\mathcal{I}C(Q;X): \;  \mbox{The space of all functions $x:Q\to X$ such that $x:Q\setminus \cup_{ j=1,\cdots,m}\{\tau_j\}\to X$ is continuous and}\\
& \qquad\qquad\qquad  \mbox{$z(\tau_j ^{+})$ and $z(\tau_j ^{-})$  exist with  $z(\tau_j )=z(\tau_j ^{-})$};\\
&P(X): \; \mbox{The set of all nonempty subsets of $X$};\\
&P_c(X): \; \mbox{The set of all closed subsets of $P(X)$};\\
&P_{k(cb)v}(X): \; \mbox{The set of all compact (closed and bounded) convex subsets of $P(X)$}.
\end{flalign*}

The norms in spaces  $C(Q;X)$, $L^{p}(Q;X)$,  and $\mathcal{I}C(Q;X)$ are respectively defined by
$$\|z\|_{C(Q;X)}=\max_{t\in Q} \|z(t)\|_{X}, \quad \|z\|_{L^{p}(Q;X)} = \left(\int_{Q} \|z(t)\|_{X} ^{p} dt \right)^{\frac{1}{p}},\;\mbox{and}\;\|z\|_{\mathcal{I}C(Q;X)}=\sup_{t\in Q} \|z(t)\|.$$

In the sequel,  let $\Gamma(\cdot)$ denote the gamma function.

\begin{definition}(\cite{FDE1})\label{sdf123}
The $q$-th fractional integral of $z(s)$ with $q>0$ is defined by\\
$$D_{0}^{-q} z(s):= \frac{1}{\Gamma(q)} \int_{0} ^{s} (s-t)^{q-1}z(t)dt,\ \ s>0.$$
\end{definition}

\begin{definition}(\cite{FDE1}) \label{d2.5}
For $\alpha\in(n-1,n)$,  the Caputo fractional order derivative of $\alpha$ of $z(s)$, denoted by $^{C}D_{0}^{\alpha}z(s)$, can  be defined by setting
$$^{C}D_{0}^{\alpha}z(s):=\frac{1}{\Gamma(n-\alpha)}\int_{0} ^{s} (s-t)^{n-\alpha-1}z^{(n)}(t)dt ,\ \ s>0.$$
\end{definition}

\begin{definition}(\cite{Carl})
The generalized directional derivative of a locally Lipschitz functional $F:Z_{2}\to R$  at $x\in Z_2$  in the direction $z\in Z_2$ and the generalized gradient of function $F$ at $v$,  denoted respectively by $F^{\circ}(x;z)$ and $\partial F(v)$,  are respectively defined by
$$F^{\circ}(x;z)=\limsup_{y\to x,\;\mu\to 0^{+}}\frac{F(y+\mu z)-F(y)}{\mu}, \quad \forall x,z\in Z_2$$
and
$$\partial F(v)=\left\{\eta\in Z_{2} ^{*}|F^{\circ}(v;z)\geq\langle \eta,z\rangle,\;\forall z\in Z_{2}\right\}, \quad \forall v\in Z_2.$$
\end{definition}

\begin{definition}(\cite{bbf2017}) \label{d2.hj5}
An operator $B:Z_{2}\to Z_{2} ^{*}$ is said to be
 \begin{itemize}
\item[(i)] monotone if $\langle Bx_1-Bx_2,x_1-x_2\rangle\geq0$, \quad $\forall x_1,x_2\in Z_2$;
\item[(ii)] strongly monotone, if there exists $m_B>0$ satisfying
$$\langle Bx_1-Bx_2,x_1-x_2\rangle\geq m_B\parallel x_1-x_2\parallel _{Z_2} ^{2},\quad  \forall x_1,x_2\in Z_2;$$
\item[(iii)] pseudomomotone, if $B$ is bounded and $x_{n}\to x$ weakly in $Z_{2}$ with
$\limsup\langle Bx_{n},x_{n}-y\rangle\leq0$ yields that
$$\liminf\langle Bx_{n},x_{n}-y\rangle\geq\langle Bx,x-y\rangle,\quad  \forall y\in Z_{2};$$
\item[(iv)] demicontinuous if $z_n\to z$ in $Z_2$ implies that $Bz_n\to Bz$ weakly in $Z_2 ^{*}$;
\item[(v)] bounded if $\Omega\subset Z_2$ is bounded implies $B(\Omega)\subset Z_2 ^{*}$ is bounded.
\end{itemize}
 \end{definition}

 \begin{definition}(\cite{bh2013}) \label{6+}
A set-valued operator $B:Z_{2}\to P(Z_{2} ^{*})$ is said to be pseudomomotone, if
 \begin{itemize}
\item[(i)] for every $x\in Z_2$, $Bx\in P_{cb}(Z_{2} ^{*})$;
\item[(ii)] for any subspace $H$ of $Z_2$, $B$ is upper semicontinuous from $H$ to $Z_{2} ^{*}$ endowed with \mbox{the weak topology}.
\item[(iii)] if $z_{n}\to z$ weakly in $Z_{2}$ and $z_n ^{*}\in Bz_n$ such that
$$\limsup\langle z_{n} ^{*},z_{n}-z\rangle\leq0,$$
then for every $x\in Z_{2}$, there exists $z ^{*}\in Bz$ such that
 $\liminf\langle z_{n} ^{*},z_{n}-x\rangle\geq\langle z ^{*},z-x\rangle$.
\end{itemize}
 \end{definition}

\begin{lemma}\label{l2wj1}
(\cite[Proposition 5.6]{han2015}) Assume that $U_1$ and $U_2$ are two reflexive Banach spaces,  $\psi:U_1\to U_2$ is a linear, continuous, and compact operator, and $\psi^{*}:U_2 ^{*}\to U_1 ^{*}$ is the adjoint operator of $\psi$. If $\varphi:U\to R$ is a locally Lipschitz functional satisfying
$$\|\partial \varphi(u)\|_{U_1 ^{*}}\leq c_{\varphi}(1+\|u\|_{U_1}),\quad \forall u\in U_1,$$
where $c_{\varphi}>0$ is a constant, then the set-valued operator $W:U_1\to P(U_1 ^{*})$, defined by
$$W(u)=\psi^{*}\partial \varphi(\psi(u)),\quad \forall u\in U_1,$$
is pseudomonotone.
 \end{lemma}

\begin{lemma}\cite[Corollary 7]{Zeng1}\label{sjxr} Assume that $U_0$  is a reflexive Banach spaces and the following conditions
 \begin{itemize}
\item[(a)] $T:U_0\to U_0 ^{*}$ is pseudomomotone and strong monotone with constant $c_T>0$;
\item[(b)] $G:U_0\to P(U_0 ^{*})$ is pseudomomotone and there exist two constants $c_G,c^{*}>0$ satisfying
$$\|G(u)\|_{U_0 ^{*}}\leq c_G\|u\|_{U_0}+c^{*}, \quad \forall u\in U_0;$$
\item[(c)]$c_G<c_T$.
\end{itemize}
hold. Then, $T+G$ is surjective in $U_0 ^{*}$, namely, for every $g\in U_0 ^{*}$, there \mbox{}exists $v\in U_0$ such that $(T+G)v\ni g$.
 \end{lemma}

Now, we rewrite   FIDHVI  as follows.
\begin{problem}\label{p1} Find $z:Q\to Z_1$ and $y:Q\to Z_2$ such that
$$
\begin{cases}
^{C}D_{0}^{\kappa}z(t)= f(t,z(t),y(t)),\quad   t\in Q,\;t\neq \tau_j,\;j=1,2,\cdots,m, \\
\Lambda z(\tau_j)=\Theta_j(z(\tau_j ^{-})),\;j=1,2,\cdots,m,\\
z(0)=z_0,\\
A(t,y(t))+N^{*}\partial J(t,Ny(t))\ni g(t,z(t)),\quad \forall t\in  Q.
\end{cases}
$$
\end{problem}

To study Problem \ref{p1}, we consider the following fractional impulsive Cauchy problem
$$
\begin{cases}
^{C}D_{0}^{\kappa}z(t)= u(t),\quad   t\in Q,\;t\neq \tau_j,\;j=1,2,\cdots,m, \\
\Delta z(\tau_j)=I_j(z(\tau_j ^{-})),\;j=1,2,\cdots,m,\\
z(0)=z_0.
\end{cases}
$$

Noting the fact that
 \begin{align*}
  z(t)= z_0-\frac{1}{\Gamma(\kappa)}\int_{0} ^{a}(a-s)^{\kappa-1}u(s)ds+\frac{1}{\Gamma(\kappa)}\int_{0} ^{t}(t-s)^{\kappa-1}u(s)ds,\quad a>0
\end{align*}
solves the Cauchy problem
$$
\begin{cases}
^{C}D_{0}^{\kappa}z(t)= u(t),\;t\in Q,\\
z(0)= z_0-\frac{1}{\Gamma(\kappa)}\int_{0} ^{a}(a-s)^{\kappa-1}u(s)ds,
\end{cases}
$$
we have the following result immediately.

\begin{lemma}\label{l2.06}
Let $\kappa\in (0,1)$ and $u\in C(Q;Z_1)$. Then the Cauchy problem
\begin{equation*}
\left\{
\begin{array}{ll}
^{C}D_{0}^{\kappa}z(t)= u(t),\;t\in Q \\
z(a)= z_0, \;a>0.
\end{array} \right.
\end{equation*}
\end{lemma}
 is equivalent to the integral equation
 \begin{align*}
  z(t)= z_0-\frac{1}{\Gamma(\kappa)}\int_{0} ^{a}(a-s)^{\kappa-1}u(s)ds+\frac{1}{\Gamma(\kappa)}\int_{0} ^{t}(t-s)^{\kappa-1}u(s)ds.
\end{align*}

Now  the following result follows.
\begin{lemma}\label{l2.07}
For $\kappa\in (0,1)$ and $u\in C(Q;Z_1)$, the Cauchy problem
\begin{equation}\label{e2.05}
\left\{
\begin{array}{ll}
^{C}D_{0}^{\kappa}z(t)= u(t),\quad   t\in Q,\;t\neq \tau_j,\;j=1,2,\cdots,m, \\
\Lambda z(\tau_j)=\Theta_j(z(\tau_j ^{-})),\;j=1,2,\cdots,m,\\
z(0)=z_0
\end{array} \right.
\end{equation}
\end{lemma}
 is equivalent to the   integral equation
 \begin{equation}\label{e2.04}
 z(t)=z_0+\sum_{i=1} ^{j}\Theta_i(z(\tau_i ^{-}))+\frac{1}{\Gamma(\kappa)}\int_{0} ^{t}(t-s)^{\kappa-1}u(s)ds,\quad \forall t\in (t_j,t_{j+1}].
\end{equation}
\textbf{Proof.} Assume that  \eqref{e2.05}  holds. If $t\in [0,\tau_1]$, then
\begin{align}\label{e2.06}
 ^{C}D_{0}^{\kappa}z(t)= u(t),\quad   t\in [0,\tau_1]\quad \mbox{with}\quad z(0)=z_0.
\end{align}
Integrating \eqref{e2.06} from $0$ to $\tau_1$, one has
\begin{align*}
 z(t)=z_0+\frac{1}{\Gamma(\kappa)}\int_{0} ^{t}(t-s)^{\kappa-1}u(s)ds.
\end{align*}
If $t\in(\tau_1,\tau_2]$, then
\begin{align}\label{e2.07}
 ^{C}D_{0}^{\kappa}z(t)= u(t),\quad   t\in (\tau_1,\tau_2]\quad \mbox{with}\quad z(\tau_1 ^{+})=z(\tau_1 ^{-})+\Theta_1(z(\tau_j ^{-}))
\end{align}
and so Lemma \ref{l2.06} implies that
\begin{align*}
 z(t)&=z(\tau_1 ^{+})-\frac{1}{\Gamma(\kappa)}\int_{0} ^{\tau_1}(\tau_1-s)^{\kappa-1}u(s)ds+\frac{1}{\Gamma(\kappa)}\int_{0} ^{t}(t-s)^{\kappa-1}u(s)ds\\
 &=z(\tau_1 ^{-})+\Theta_1(z(\tau_1 ^{-}))-\frac{1}{\Gamma(\kappa)}\int_{0} ^{\tau_1}(\tau_1-s)^{\kappa-1}u(s)ds+\frac{1}{\Gamma(\kappa)}\int_{0} ^{t}(t-s)^{\kappa-1}u(s)ds\\
 &=z_0+\Theta_1(z(\tau_1 ^{-}))+\frac{1}{\Gamma(\kappa)}\int_{0} ^{t}(t-s)^{\kappa-1}u(s)ds.
\end{align*}
If $t\in(\tau_2,\tau_3]$, then using Lemma \ref{l2.06} again,  we have
\begin{align*}
 z(t)&=z(\tau_2 ^{+})-\frac{1}{\Gamma(\kappa)}\int_{0} ^{\tau_2}(\tau_1-s)^{\kappa-1}u(s)ds+\frac{1}{\Gamma(\kappa)}\int_{0} ^{t}(t-s)^{\kappa-1}u(s)ds\\
 &=z(\tau_2 ^{-})+\Theta_2(z(\tau_2 ^{-}))-\frac{1}{\Gamma(\kappa)}\int_{0} ^{\tau_2}(\tau_1-s)^{\kappa-1}u(s)ds+\frac{1}{\Gamma(\kappa)}\int_{0} ^{t}(t-s)^{\kappa-1}u(s)ds\\
 &=z_0+\Theta_1(z(\tau_1 ^{-}))+\Theta_2(z(\tau_2 ^{-}))+\frac{1}{\Gamma(\kappa)}\int_{0} ^{t}(t-s)^{\kappa-1}u(s)ds.
\end{align*}
Similarly, if $t\in(\tau_j,\tau_{j+1}]$, then we can show that
\begin{equation*}
 z(t)=z_0+\sum_{i=1} ^{j}\Theta_i(z(\tau_i ^{-}))+\frac{1}{\Gamma(\kappa)}\int_{0} ^{t}(t-s)^{\kappa-1}u(s)ds,\quad \forall  t\in (t_j,t_{j+1}].
\end{equation*}

Conversely, suppose that  \eqref{e2.04} holds. If $t\in (0,\tau_1]$, then we know that \eqref{e2.05} holds by the fact that $^{C}D_{0}^{\kappa}$ is the  inverse of $D_{0}^{-\kappa}$.  If $t\in(\tau_j,\tau_{j+1}], j= 1,2,...,m$, since the Caputo fractional derivative for a constant is  zero,  one has $^{C}D_{0}^{\kappa}z(t)= u(t),t\in(\tau_j,\tau_{j+1}]$ and $\Lambda z(\tau_j)=\Theta_j(z(\tau_j ^{-}))$.  \hfill$\Box$

From Lemma \ref{l2.07}, we have the following definition.
\begin{definition}\label{sdf123}
 A pair $(z,y)\in \mathcal{I}C(Q;Z_{1}) \times \mathcal{I}C(Q;Z_2)$  is said to be a solution of Problem \ref{p1} if it satisfies  the following system
\begin{align}
&z(t)=z_0+\sum_{i=1} ^{j}\Theta_i(z(\tau_i ^{-}))+\frac{1}{\Gamma(\kappa)}\int_{0} ^{t}(t-s)^{\kappa-1}f(s,z(s),y(s))ds,\quad \forall t\in (t_j,t_{j+1}],j= 1,2,...,m,\label{bh2.2}\\
&A(t,y(t))+\eta\ni g(t,z(t)), \quad \forall t\in  Q,\label{bh2.3-}\\
&\eta\in N^{*}\partial J(t,Ny(t)),\quad \forall t\in  Q.\label{bh2.c4-}
\end{align}
\end{definition}

Finally, we recall the following nonlinear impulsive Gronwall inequality.

\begin{lemma}\cite[Lemma 3.4]{imgronwall}\label{la2.77aa}
Let $z\in \mathcal{I}C(Q;Z_1)$ satisfy the \mbox{following inequality}
\begin{align*}
  \|z(t)\|\leq k_1+k_2\int_{0} ^{t}(t-s)^{\kappa-1}\|z(s)\|ds+\sum_{0<\tau_{j}<t}d_j\|z(\tau_j ^{-})\|,
\end{align*}
where $ k_1,k_2, d_j\geq 0$ are constants. Then
\begin{align*}
  \|z(t)\|\leq k_1\left[1+D^{*}E_{\kappa}(k_2\Gamma(\kappa)t^{\kappa})\right]^{j}E_{\kappa}(k_2\Gamma(\kappa)t^{\kappa}),\quad \forall t\in (t_j,t_{j+1}],
\end{align*}
where $D^{*}=\max\{d_j:j=1,\cdots,m\}$ and $E_{\gamma}$ is the \mbox{Mittag-Leffler} function \cite{FDE1} defined by
\begin{align*}
  E_{\gamma}(h)=\sum_{j=0} ^{\infty}\frac{h^{j}}{\Gamma(\gamma h+1)}, \quad h\in \mathbb{C}, \quad \mbox{Re}(\gamma) >0.
\end{align*}
 \end{lemma}

\section{Existence and uniqueness}\label{sec.3}
\setcounter{equation}{0}
To study the solvability of Problem \ref{p1}, we need the following assumptions.

\begin{itemize}
\item[(Hf)]  $f:Q\times Z_{1}\times Z_{2}\to Z_{1}$ is a map such that
\begin{itemize}
\item[(i)] for any given $(z,y)\in Z_{1}\times Z_{2}$, $f(\cdot,z,y)$ is continuous;
\item[(ii)] there exists  $M_1>0$ satisfying
$$\|f(t,z_{1},y_{1})-f(t,z_{2},y_{2})\|_{Z_{1}}\leq M_1\left(\|z_{1}-z_{2}\|_{ Z_{1}}+\|y_{1}-y_{2}\|_{ Z_{2}}\right), \; \forall (t,z_{i},y_{i})\in Q \times Z_{1}\times Z_{2}, \; i=1,2;$$
\item[(iii)] there exists  $\phi\in  L_{+} ^{\frac{1}{p}}[0,T](0<p<\kappa<1)$ satisfying
$$\|f(t,z,y)\|_{Z_{1}}\leq \phi(t), \quad \forall (t,z,y)\in Q \times Z_{1}\times Z_{2}.$$
\end{itemize}

\item[(HI)] For each $j\in \{1,2,...,m\}$,  $\Theta_j:Z_1\to Z_1$ is bounded and there exists  $d_j> 0$ satisfying
$$\|\Theta_j(z_1)-\Theta_j(z_1)\|_{Z_1}\leq d_j\|z_1-z_2\|_{Z_1},\quad \forall z_1,z_2\in Z_1.$$

\item[(HA)]  $A:Q\times Z_2\to Z_{2} ^{*}$  is a map such that
\begin{itemize}
\item[(i)]for any given $y\in Z_2$, $A(\cdot,y)$ is continuous;
\item[(ii)]for any given $t\in Q$, $A(t,\cdot)$ is bounded, demicontinuous, and strongly monotone with the constant $m_A$.
\end{itemize}

\item[(HN)] $N\in L(Z_{2},Y)$ is a compact and surjective operator.

\item[(HJ)] $J:Q\times Y\to R$ is a functional satisfying
\begin{itemize}
\item[(i)] for any given $x\in Y$, $J(\cdot,x)$ is continuous;
\item[(ii)] for any given $t\in Q$, $J(t,\cdot)$ is locally Lipschitz;
\item[(iii)] there exists  $m_{J}>0$ satisfying
$$ \|\partial J(t,x)\|_{Z_2 ^{*}}\leq m_{J}(\|x\|_{Y}+1), \forall (t,x)\in Q\times Y.$$
\item[(iv)] there exists  $c_{J}>0$ satisfying
$$\langle \theta_1-\theta_2, y_1-y_2\rangle\geq -c_{J}\|y_1-y_2\|_{Y} ^{2},\quad \forall \theta_i\in \partial J(t,y_i), \;(t,y_i)\in Q\times Y, i=1,2.$$
\end{itemize}

\item[(Hg)] $g:Q\times Z_1\to Z_{2} ^{*}$  is a map such that
\begin{itemize}
\item[(i)]for any given $z\in Z_{1}$, $g(\cdot,z)$ is continuous;
\item[(ii)] there exists  $m_g>0$ satisfying
$$\|g(t,z_{1})-g(t,z_{2})\|_{Z_{2} ^{*}}\leq m_g\|z_{1}-z_{2}\|_{ Z_{1}}, \quad \forall (t,z_i)\in Q \times Z_{1}, i=1,2.$$
\end{itemize}

\item[(HO)]
\begin{itemize}
\item[(i)] $m_A> c_{J}\|N\|^{2}$, where $\|N\|=\|N\|_{L(Z_{2},Y)}$;
\item[(ii)] $\frac{T^{\kappa}M_1 m_g}{\kappa(m_A- c_{J}\|N\|^{2})\Gamma(\kappa)}<1$.
\end{itemize}
\end{itemize}

We first consider nonlinear inclusions \eqref{bh2.3-}-\eqref{bh2.c4-}.

\begin{lemma}\label{l3.3rf}
For any given $z\in \mathcal{I}C(Q;Z_{1})$, nonlinear inclusions \eqref{bh2.3-}-\eqref{bh2.c4-} have a unique solution $y\in \mathcal{I}C(Q;Z_2)$ providing that assumptions (HA),(HN), (HJ), (Hg), and (HO) hold. Moreover, for any $z_{1},z_{2}\in \mathcal{I}C(Q;Z_{1})$,  one has
\begin{equation}\label{e1.sd2}
 \|y_1(t)-y_2(t)\|_{Z_2}\leq \frac{m_g}{m_A- c_{J}\|N\|^{2}}\|z_1(t)-z_2(t)\|_{Z_1},\quad \forall  t\in Q,
\end{equation}
where $y_{1},y_{2}\in \mathcal{I}C(Q;Z_2)$ are the solutions of \eqref{bh2.3-}-\eqref{bh2.c4-} with respect to $z_1$ and $z_2$, respectively.
\end{lemma}
\textbf{Proof.} For given $z\in \mathcal{I}C(Q;Z_{1})$ and  $t\in Q$, define two operators $\widehat{A}:Z_2\to Z_2 ^{*}$ and $\widehat{N}:Z_2\to P(Z_2 ^{*})$ as follows:
$$\widehat{A}y=A(t,y), \quad \widehat{N}y=N^{*}\partial J^{\circ}(t,,Ny),\qquad \forall (t,y)\in Q\times Z_2.$$
 For simplicity, we do \mbox{not indicate their dependence} $t$. Using  (HA), (HN), (HJ),  (HO),  Lemma \ref{l2wj1} and \cite[Lemma 3]{bbf2017}, we deduce that the operators $\widehat{A}$ and $\widehat{N}$  are  pseudomomotone and
 \begin{align*}
  \|\widehat{N}x\|_{E_2 ^{*}}& \leq  \|N^{*}\|\| \partial J^{\circ}(t,Nx)\|\\
  &\leq \|N^{*}\|\left(m_{J_1}\|Nx\|_{X}+m_{J_2}\right)\\
   &\leq m_{J_1}\|N\|^{2}\|x\|_{Z_2}+m_{J_2}\|N\|,\quad \forall x\in Z_2.
 \end{align*}
 By applying  Lemma \ref{sjxr} with $B=\widehat{A}$ and $A=\widehat{N}$, we know that inclusions \eqref{bh2.3-}-\eqref{bh2.c4-} has a solution $y(t)$ for all $t\in Q$.
 Next we show that the solution $y(t)$ is unique. Let $y_1,y_2\in Z_2$ be solutions to  \eqref{bh2.3-}-\eqref{bh2.c4-}. Then, there exist $\eta_1,\eta_2\in N^{*}\partial J^{\circ}(t,Ny_i(t))$   satisfying
 \begin{align*}
   A(t,y_i)+\eta_i= g(t,z(t)),\quad  \quad \forall  t\in Q, \; i=1,2.
 \end{align*}
 Subtracting the above two equations and taking the result in duality with $y_1-y_2$, we have
 \begin{align*}
   \langle A(t,y_1)-A(t,y_1),y_1-y_2\rangle_{ Z_2 ^{*}\times Z_2}=\langle \eta_2-\eta_1,y_1-y_2\rangle_{ Z_2 ^{*}\times Z_2}.
 \end{align*}
By assumptions (HA) and (HJ), one has $\left(m_A- c_{J}\|N\|^{2}\right)\|y_1-y_2\|_{Z_2} ^{2}\leq 0$
and so the assumption (HO) implies that $y_1=y_2$, which is our claim.

 In what follows, we start  by  showing that \eqref{e1.sd2} holds.  Let $z_i(t)\in Z_1(i=1,2)$ and denote $z_{i}(t)=z_{i}$, $y_{i}(t)=y_{i}$, $g(t,z_{i}(t))=g_{i}$  with $i=1,2$. It follows from  \eqref{bh2.3-} and \eqref{bh2.c4-} that
\begin{align}
 A(t,y_1)+\varsigma_1= g_1,\quad \varsigma_1\in N^{*}\partial J^{\circ}(t,,Ny_1)\label{a3.3a2}
\end{align}
and
\begin{align}
A(t,y_2)+\varsigma_2= g_2,\quad \varsigma_2\in N^{*}\partial J^{\circ}(t,Ny_2).\label{a3.4a2}
 \end{align}
  Subtracting \eqref{a3.4a2} from \eqref{a3.3a2} and taking the result in duality with $y_1-y_2$, we have
  \begin{align*}
   \langle A(t,y_1)-A(t,y_1),y_1-y_2\rangle_{ Z_2 ^{*}\times Z_2}+\langle \varsigma_1-\varsigma_2,y_1-y_2\rangle_{ Z_2 ^{*}\times Z_2}=\langle g_1-g_2,y_1-y_2 \rangle_{ Z_2 ^{*}\times Z_2}.
 \end{align*}
 By assumptions (HJ), (HA) and (Hg), one has
 \begin{align*}
   \left(m_A- c_{J}\|N\|^{2}\right)\|y_1-y_2\|_{Z_2} ^{2}\leq \|g_1-g_2\|_{Z_2 ^{*}} \|y_1-y_2\|_{Z_2}\leq m_g\|z_1-z_2\|_{Z_1}\|y_1-y_2\|_{Z_2}.
 \end{align*}
 Thus, the assumption (HO) implies that
  \begin{align}\label{e3.04a}
   \|y_1-y_2\|_{Z_2}\leq \frac{m_g}{m_A- c_{J}\|N\|^{2}}\|z_1-z_2\|_{Z_1}.
 \end{align}
 It follows from \eqref{e3.04a} that the map $Z_1\ni z(t)\mapsto y(t)\in Z_2$ is continuous for all $t\in Q$. Since $z\in \mathcal{I}C(Q;Z_{1})$, we know that  $y\in  \mathcal{I}C(Q;Z_2)$. By \eqref{e3.04a}, we conclude that,  for any given $z\in \mathcal{I}C(Q;Z_{1})$,  nonlinear inclusions \eqref{bh2.3-}-\eqref{bh2.c4-} has a unique solution $y\in \mathcal{I}C(Q;Z_2)$. Moreover, for any given $z_{1},z_{2}\in \mathcal{I}C(Q;Z_{1})$,  \eqref{e1.sd2} holds due to \eqref{e3.04a}.  \hfill$\Box$

\begin{theorem}\label{hv+zxdl}
 Problem \ref{p1} admits a unique solution $(z,y)\in \mathcal{I}C(Q;Z_{1}) \times \mathcal{I}C(Q;Z_{2})$ providing that assumptions (HA), (Hf), (HI), (HN), (HJ), (Hg), and (HO) hold.
 \end{theorem}

\textbf{Proof.} For any given $z\in \mathcal{I}C(Q;Z_{1})$, Lemma \ref{l3.3rf} shows that nonlinear inclusions \eqref{bh2.3-}-\eqref{bh2.c4-} admits a unique solution $y_{z}$.  Define an operator $\Sigma:\mathcal{I}C(Q;Z_{1})\to \mathcal{I}C(Q;Z_{1})$ by setting
$$\Sigma z(t)=z_0+\sum_{i=1} ^{j}\Theta_i(z(\tau_i ^{-}))+\frac{1}{\Gamma(\kappa)}\int_{0} ^{t}(t-s)^{\kappa-1}f(s,z(s),y_z(s))ds.$$
Then the assumption (Hf) implies that $\Sigma$ is well defined. To prove Theorem \ref{hv+zxdl}, we only need to show that  $\Sigma$ admits a unique fixed point in $\mathcal{I}C(Q;Z_{1})$.

To this end, we first show that $\Sigma z\in \mathcal{I}C(Q;Z_{1})$ for any $z\in \mathcal{I}C(Q;Z_{1})$.  In fact, let $z\in C([0,\tau_1],Z_1)$ and $\iota>0$ be given. When $t\in [0,\tau_1]$, by the H\"{o}lder inequality and the assumption (Hf), we have
\begin{align*}
&\quad \; \|(\Sigma z)(t+\iota)-(\Sigma z)(t)\|_{Z_{1}}\\
&\leq \frac{1}{\Gamma(\kappa)}\int_{0} ^{t}\left((t-s)^{\kappa-1}-(t+\iota-s)^{\kappa-1}\right)\|f(s,z(s),y_{z}(s)\|_{Z_{1}}ds\\
&\quad +\frac{1}{\Gamma(\kappa)}\int_{t} ^{t+\iota}(t+\iota-s)^{\kappa-1}\|f(s,z(s),y_{z}(s)\|_{Z_{1}}ds\\
&\leq \frac{1}{\Gamma(\kappa)}\int_{0} ^{t}\left((t-s)^{\kappa-1}-(t+\iota-s)^{\kappa-1}\right)\phi(s)ds+\frac{1}{\Gamma(\kappa)}\int_{t} ^{t+\iota}(t+\iota-s)^{\beta-1}\phi(s)ds\\
&\leq \frac{1}{\Gamma(\kappa)}\left(\int_{0} ^{t}\left((t-s)^{\kappa-1}-(t+\iota-s)^{\kappa-1}ds\right)^{\frac{1}{1-p}}\right)^{1-p}\left(\int_{0} ^{t}\left(\phi(s)\right)^{p}ds\right)^{\frac{1}{p}}\\
&\quad +\frac{1}{\Gamma(\kappa)}\left(\int_{t} ^{t+\iota}\left((t+\iota-s)^{\kappa-1}\right)^{\frac{1}{1-p}}ds\right)^{1-p}\left(\int_{t} ^{t+\iota}\left(\phi(s)\right)^{\frac{1}{p}}ds\right)^{p}\\
&\leq \frac{M}{\Gamma(\kappa)}\left(\int_{0} ^{t}\left((t-s)^{\kappa}-(t+\iota-s)^{\alpha}\right)ds\right)^{1-p}+\frac{M}{\Gamma(\kappa)}\left(\int_{t} ^{t+\iota}(t+\iota-s)^{\alpha}ds\right)^{1-p}\\
&\leq  \frac{M}{\Gamma(\kappa)(1+\alpha)^{1-p}}\left(|(t+\iota)^{1+\alpha}-t^{1+\alpha}|+\iota^{1+\alpha}\right)^{1-p}
+\frac{M}{\Gamma(\kappa)(1+\alpha)^{1-p}}\iota^{(1+\alpha)(1-p)}\\
&\leq \frac{2M}{\Gamma(\kappa)(1+\alpha)^{1-p}}\iota^{(1+\alpha)(1-p)}+\frac{M}{\Gamma(\kappa)(1+\alpha)^{1-p}}\iota^{(1+\alpha)(1-p)}\\
&\leq\frac{3M}{\Gamma(\kappa)(1+\alpha)^{1-p}}\iota^{(1+\alpha)(1-p)}\\
& \to 0
 \end{align*}
as $\iota\to 0$, where $M=\|\phi\|_{L^{\frac{1}{p}}[0,T]}$ and $\alpha=\frac{\kappa-1}{1-p}\in(-1,0)$. This shows that $\Sigma z\in C([0,\tau_1],Z_1)$.

When $t\in(\tau_1,\tau_2]$, using the same argument, one has
\begin{align*}
  \|(\Sigma z)(t+\iota)-(\Sigma z)(t)\|_{Z_{1}}\leq\frac{3M}{\Gamma(\kappa)(1+\alpha)^{1-p}}\iota^{(1+\alpha)(1-p)} \to 0 \quad \mbox{as} \; \iota\to 0,
\end{align*}
which implies that $\Sigma z\in C((\tau_1,\tau_2],Z_1)$.

Similarly, when $t\in(\tau_j,\tau_{j+1}],\;j=1,2,\cdots,m$, we can show that
\begin{align*}
  \|(\Sigma z)(t+\iota)-(\Sigma z)(t)\|_{Z_{1}}\leq\frac{3M}{\Gamma(\kappa)(1+\alpha)^{1-p}}\iota^{(1+\alpha)(1-p)} \to 0 \quad \mbox{as} \; \iota\to 0
\end{align*}
and so $\Sigma z\in C((\tau_j,\tau_{j+1}],Z_1)$.

Combining all the above  we see that  $\Sigma z\in \mathcal{I}C(Q;Z_{1})$ for any $z\in \mathcal{I}C(Q;Z_{1})$.

Next we prove that $\Sigma$ is a contractive map.  For given $z_1,z_2\in \mathcal{I}C(Q;Z_{1})$, by the assumption (Hf), it follows from (\ref{e1.sd2}) that
\begin{align*}
&\quad \|(\Sigma z_1)(t)-(\Sigma z_2)(t)\|_{Z_{1}}\\
&\leq \frac{1}{\Gamma(\kappa)}\int_{0} ^{t}(t-s)^{\kappa-1}\|f(s,z_1(s),y_{z_1}(s)- f(s,z(s),y_{z_1}(s)\|_{Z_{1}}ds\\
&\leq \frac{M_1}{\Gamma(\kappa)}\int_{0} ^{t}(t-s)^{\kappa-1}\left(\|z_1(s)-z_2(s)\|_{Z_{1}}+\|y_{z_1}(s)-y_{z_2}(s)\|_{Z_{2}}\right)ds\\
&\leq \frac{M_1 m_g}{(m_A- c_{J}\|N\|^{2})\Gamma(\kappa)}\int_{0} ^{t}(t-s)^{\kappa-1}\|z_1(s)-z_2(s)\|_{Z_{1}}ds\\
&\leq \frac{M_1 m_g}{(m_A- c_{J}\|N\|^{2})\Gamma(\kappa)}\left(\int_{0} ^{t}(t-s)^{\kappa-1}ds\right)\|z_1-z_2\|_{\mathcal{I}C(Q;Z_{1})}\\
&\leq \frac{T^{\beta}M_1 m_g}{\beta(m_A- c_{J}\|N\|^{2})\Gamma(\kappa)}\|z_1-z_2\|_{\mathcal{I}C(Q;Z_{1})}
 \end{align*}
and so
\begin{align*}
\|\Sigma z_1-\Sigma z_2\|_{\mathcal{I}C(Q;Z_{1})}\leq \frac{T^{\beta}M_1 m_g}{\beta(m_A- c_{J}\|N\|^{2})\Gamma(\kappa)}\|z_1-z_2\|_{\mathcal{I}C(Q;Z_{1})}.
\end{align*}
Now the assumption (HO) implies that $\Sigma$ is a contractive map,  and so $\Sigma$ admits a unique solution $z\in\mathcal{I}C(Q;Z_{1})$ by employing the Banach fixed point theorem. \hfill$\Box$

\section{ A convergence result}\label{sec.4}

We investigate the perturbation problem of Problem \ref{p1}  to prove a convergence result which describes the stability of the solution in relation to perturbation data. To this end, let $\delta>0$ and $J_{\delta}$ be the perturbed data of $J$ such that $J_{\delta}$ satisfies assumptions (HJ) and (HO). More precisely, we examine the following perturbation problem of Problem \ref{p1}: find a pair of functions $(z_{\delta},y_{\delta})\in \mathcal{I}C(Q;Z_{1}) \times \mathcal{I}C(Q;Z_{2})$  such that
\begin{eqnarray}\label{e4.05}
\left\{
\begin{array}{l}
^{C}D_{0}^{\kappa}z_{\delta}(t)= f(t,z_{\delta}(t),y_{\delta}(t)),\quad   t\in Q,\;t\neq \tau_j,\;j=1,2,\cdots,m, \\
\Lambda z_{\delta}(\tau_j)=\Theta_j(z_{\delta}(\tau_j ^{-})),\;j=1,2,\cdots,m,\\
z(0)=z_0,\\
\langle A(t,y_{\delta}(t)),x\rangle+J_{\delta} ^{\circ}(t,y_{\delta}(t),Ny_{\delta}(t);Nx)\geq\langle g(t,z_{\delta}(t)),x\rangle,\quad \forall (t,x)\in  Q\times Z_2.
\end{array}
\right.
\end{eqnarray}

We denote the constants involved in the assumption (HJ) by $m_{J\delta}$ and $c_{J_{\delta}}$. Furthermore, we introduce the following  assumptions.
\begin{itemize}
\item[(HJ$^{*}$)] $J_{\delta}:Q\times Y\to R$ is a functional satisfying
\begin{itemize}
\item[(i)]  there exists a function $V:R^{+}\to R^{+}$ satisfying, for any $(t,y)\in Q\times Z_2$ and $\delta>0$,
$$\|\zeta-\zeta_{\delta}\|_{Z_2 ^{*}}\leq V(\delta),\quad \forall (\zeta,\zeta_{\delta})\in N^{*}\partial J(t,Ny(t))\times N^{*}\partial J_{\delta}(t,Ny(t));$$
\item[(ii)]$\lim_{\delta\to 0} V(\delta)=0.$
\end{itemize}

\item[(HO$^{*}$)] There exists  $m_{A0}>0$ such that
\begin{itemize}
\item[(i)] $m_A>m_{A0}> c_{J\delta}\|N\|^{2}$, where $\|N\|=\|N\|_{L(Z_{2},Y)}$;
\item[(ii)] $\frac{T^{\kappa}M_1 m_g}{\kappa(m_A- c_{J\delta}\|N\|^{2})\Gamma(\kappa)}<1.$
\end{itemize}
\end{itemize}

\begin{theorem}\label{hvpre5dl}
Suppose that assumptions (HA), (Hf), (HI), (HN), (HJ), (Hg), (HO), (HO$^{*}$) and (HJ$^{*}$) hold. Then
\begin{itemize}
\item[(i)] for each $\delta>0$, the perturbation problem \eqref{e4.05} has a unique solution $(z_{\delta},y_{\delta})\in \mathcal{I}C(Q;Z_{1}) \times \mathcal{I}C(Q;Z_{2})$;
\item[(ii)]$(z_{\delta},y_{\delta})$ converges to $(z(t),y(t))$, the solution  of Problem \ref{p1}, i.e.,
\begin{align}\label{e4900}
(z_{\delta}(t),y_{\delta}(t))\to (z(t),y(t))\quad \mbox{as}\;\delta\to 0, \;\forall t\in Q.
\end{align}
\end{itemize}
 \end{theorem}
\textbf{Proof.} (i) In view of Theorem \ref{hv+zxdl}, the proof is obvious.

(ii) By Definition \ref{sdf123}, we consider the  problem:
\begin{align}
&z_{\delta}(t)=z_0+\sum_{i=1} ^{j}\Theta_i(z_{\delta}(\tau_i ^{-}))+\frac{1}{\Gamma(\kappa)}\int_{0} ^{t}(t-s)^{\kappa-1}f(s,z_{\delta}(s),y_{\delta}(s))ds,\quad \forall t\in (t_j,t_{j+1}],j= 1,2,...,m,\label{bch2.2}\\
&A(t,y_{\delta}(t))+\eta_{\delta}\ni g(t,z_{\delta}(t)),\label{bh2v.3-}\\
&\eta_{\delta}\in N^{*}\partial J_{\delta}(t,Ny_{\delta}(t)),\quad \forall (t,x)\in  Q\times Z_2.\label{bh2.4-}
\end{align}
Subtracting \eqref{bh2v.3-} from \eqref{bh2.3-} and multiplying the result by $y(t)-y_{\delta}(t)$, we have
\begin{align*}
  &\quad  \langle A(t,y(t))-A(t,y_{\delta}(t)),y(t)-y_{\delta}(t)\rangle_{ Z_2 ^{*}\times Z_2}+\langle\eta-\eta_{\delta},y(t)-y_{\delta}(t)\rangle_{ Z_2 ^{*}\times Z_2}\\
  &=\langle g(t,z(t))-g(t,z_{\delta}(t)),y(t)-y_{\delta}(t) \rangle_{ Z_2 ^{*}\times Z_2},\quad  \forall (t,\eta,\eta_{\delta})\in Q\times N^{*}\partial J(t,Ny(t))\times N^{*}\partial J_{\delta}(t,Ny_{\delta}(t)).
 \end{align*}
 Since
 \begin{align*}
  \langle\eta-\eta_{\delta},y(t)-y_{\delta}(t)\rangle_{ Z_2 ^{*}\times Z_2}=&\langle\eta-\xi_{\delta},y(t)-y_{\delta}(t)\rangle_{ Z_2 ^{*}\times Z_2}+\langle\xi_{\delta}-\eta_{\delta},y(t)-y_{\delta}(t)\rangle_{ Z_2 ^{*}\times Z_2}, \\
   & \forall (t,\eta,\xi_{\delta},\eta_{\delta})\in Q\times N^{*}\partial J(t,Ny(t))\times  N^{*}\partial J(t,Ny_{\delta}(t))\times N^{*}\partial J_{\delta}(t,Ny_{\delta}(t)),
 \end{align*}
 one has
 \begin{align*}
  &\quad  \langle A(t,y(t))-A(t,y_{\delta}(t)),y(t)-y_{\delta}(t)\rangle_{ Z_2 ^{*}\times Z_2}+\langle\eta-\xi_{\delta},y(t)-y_{\delta}(t)\rangle_{ Z_2 ^{*}\times Z_2}\\
  &= \langle\eta_{\delta}-\xi_{\delta},y(t)-y_{\delta}(t)\rangle_{ Z_2 ^{*}\times Z_2}+ \langle g(t,z(t))-g(t,z_{\delta}(t)),y(t)-y_{\delta}(t) \rangle_{ Z_2 ^{*}\times Z_2},\\
    &\qquad \forall (t,\eta,\xi_{\delta},\eta_{\delta})\in Q\times N^{*}\partial J(t,Ny(t))\times  N^{*}\partial J(t,Ny_{\delta}(t))\times N^{*}\partial J_{\delta}(t,Ny_{\delta}(t)).
 \end{align*}
Note that the assumption (HA) implies that
\begin{align}\label{e4.013}
  \langle A(t,y(t))-A(t,y_{\delta}(t)),y(t)-y_{\delta}(t)\rangle_{ Z_2 ^{*}\times Z_2}\geq m_A\|y(t)-y_{\delta}(t)\|_{Z_2} ^{2}\quad \forall t\in Q.
\end{align}
Using  assumptions (HJ) and (HJ$^{*}$), we obtain
\begin{align}\label{e4.015}
 &\quad\langle\eta-\xi_{\delta},y(t)-y_{\delta}(t)\rangle_{ Z_2 ^{*}\times Z_2}\notag\\
 &\geq -c_{J}\|N\|^{2}\|y(t)-y_{\delta}(t)\|_{Z_2} ^{2},\;
 \forall (t,\eta,\xi_{\delta})\in Q\times N^{*}\partial J(t,Ny(t))\times  N^{*}\partial J(t,Ny_{\delta}(t))
\end{align}
and
\begin{align}\label{e4.014}
 &\quad\langle\eta_{\delta}-\xi_{\delta},y(t)-y_{\delta}(t)\rangle_{ Z_2 ^{*}\times Z_2}\notag\\
 &\leq \|\eta_{\delta}-\xi_{\delta}\|_{Z_2 ^{*}}\|y(t)-y_{\delta}(t)\|_{Z_2}\notag\\
 &\leq V(\delta)\|y(t)-y_{\delta}(t)\|_{Z_2},
 \; \forall (t,\xi_{\delta},\eta_{\delta})\in Q\times  N^{*}\partial J(t,Ny_{\delta}(t))\times N^{*}\partial J_{\delta}(t,Ny_{\delta}(t)).
\end{align}
We conclude from the assumption (Hg) that
\begin{align}\label{e4.016}
 &\quad \langle g(t,z(t))-g(t,z_{\delta}(t)),y(t)-y_{\delta}(t) \rangle_{ Z_2 ^{*}\times Z_2}\notag\\
 &\leq \| g(t,z(t))-g(t,z_{\delta}(t))\|_{Z_2 ^{*}}\|y(t)-y_{\delta}(t)\|_{Z_2}\notag\\
 &\leq m_g\|y(t)-y_{\delta}(t)\|_{Z_2}\|z(t)-z_{\delta}(t)\|_{Z_1} ,\quad  \forall t\in Q.
\end{align}
Combining \eqref{e4.013}-\eqref{e4.016}, one has
\begin{align*}
 \left(m_A- c_{J}\|N\|^{2}\right)\|y(t)-y_{\delta}(t)\|_{Z_2} ^{2}\leq V(\delta)\|y(t)-y_{\delta}(t)\|_{Z_2} +m_g\|y(t)-y_{\delta}(t)\|_{Z_2}\|z(t)-z_{\delta}(t)\|_{Z_1}.
\end{align*}
Thus, the assumption (HO)  yields that
\begin{align}\label{e4.017}
 \|y(t)-y_{\delta}(t)\|_{Z_2} \leq \frac{V(\delta)}{m_A- c_{J}\|N\|^{2}} +\frac{m_g}{m_A- c_{J}\|N\|^{2}}\|z(t)-z_{\delta}(t)\|_{Z_1}.
\end{align}
Subtracting \eqref{bch2.2} from \eqref{bh2.2}, by assumptions (Hf), (HI) and estimation \eqref{e4.017}, one has
\begin{align*}
&\quad \|z_{\delta}(t)-z(t)\|_{Z_{1}}\\
&\leq \frac{1}{\Gamma(\kappa)}\int_{0} ^{t}(t-s)^{\kappa-1}\|f(s,z_{\delta}(s),y_{\delta}(s)- f(s,z(s),y(s)\|_{Z_{1}}ds+\sum_{i=1} ^{j}\|\Theta_i(z_{\delta}(\tau_i ^{-}))-\Theta_i(z(\tau_i ^{-}))\|_{Z_{1}}\\
&\leq \frac{M_1}{\Gamma(\kappa)}\int_{0} ^{t}(t-s)^{\kappa-1}\left(\|z(t)-z_{\delta}(t)\|_{Z_1}+\|y(t)-y_{\delta}(t)\|_{Z_2}\right)ds+
\sum_{i=1} ^{j}d_j\|z_{\delta}(\tau_i ^{-})-z(\tau_i ^{-})\|_{Z_{1}}\\
&\leq  \frac{M_1}{\Gamma(\kappa)}\int_{0} ^{t}(t-s)^{\kappa-1}
\left[\frac{V(\delta)}{m_A- c_{J}\|N\|^{2}} +\left(\frac{m_g}{m_A- c_{J}\|N\|^{2}}+1\right)\|z(t)-z_{\delta}(t)\|_{Z_1}\right]ds\\
&\quad +\sum_{i=1} ^{j}d_j\|z_{\delta}(\tau_i ^{-})-z(\tau_i ^{-})\|_{Z_{1}}\\
&\leq  \frac{T^{\kappa}M_1}{\kappa\Gamma(\kappa)(m_A- c_{J}\|N\|^{2})}V(\delta)+\frac{M_1}{\Gamma(\kappa)}\left(\frac{m_g}{m_A- c_{J}\|N\|^{2}}+1\right)
\int_{0} ^{t}(t-s)^{\kappa-1}\|z(t)-z_{\delta}(t)\|_{Z_1}ds\\
&\quad +\sum_{i=1} ^{j}d_j\|z_{\delta}(\tau_i ^{-})-z(\tau_i ^{-})\|_{Z_{1}}.
\end{align*}
Now by Lemma \ref{la2.77aa} with $k_1= \frac{T^{\kappa}M_1}{\kappa\Gamma(\kappa)(m_A- c_{J}\|N\|^{2})}V(\delta)$ and $k_2=\frac{M_1}{\Gamma(\kappa)}\left(\frac{m_g}{m_A- c_{J}\|N\|^{2}}+1\right)$, there exists $H^{*}>0$ such that
\begin{align*}
 \|z_{\lambda}(t)-z(t)\|_{Z_{1}}\leq H^{*}V(\delta),
\end{align*}
where $H^{*}$ is independent of $z, z_{\lambda},y, y_{\lambda}$ and $t$.
By the assumption (HJ$^{*}$), we assert that
\begin{align}\label{e4.018}
 \|z_{\lambda}(t)-z(t)\|_{Z_{1}}\to 0 \; \mbox{as}\; \delta\to 0.
\end{align}
Using \eqref{e4.017}, \eqref{e4.018} and  the assumption (HJ$^{*}$), one has
\begin{align}\label{e4.019}
 \|y(t)-y_{\delta}(t)\|_{Z_2}\to 0 \; \mbox{as}\; \delta\to 0.
\end{align}
Thus, (\ref{e4.018}) and (\ref{e4.019}) finish the proof.
 \hfill$\Box$

\section{ An application}\label{sec.5}

In this section we show that the results obtained in Sections 3 and 4 can be applied to study  the frictional contact problem (Problem \ref{p2}) between an elastic body and a foundation over time interval $Q$. We suppose that the surface traction may change suddenly in a short time, such as shocks, and consequently, which can be described by a fractional impulsive differential equations. We show that the weak form of  Problem \ref{p2} leads to  Problem \ref{p1} analyzed in Sections 3 and 4. Then Theorems \ref{hv+zxdl} and  \ref{hvpre5dl} are applied to obtain the unique solvability of the frictional contact problem mentioned above as well as the convergence result of the perturbation problem.

We shortly review the basic notations and its mechanical interpretations. A deformable elastic body  occupies a regular Lipschitz domain  $V\in R^{n} (n=2,3)$ with the boundary $\partial V$.  The boundary $\partial V$ consists of three measurable disjoint parts $\Sigma_{1}$, $\Sigma_{2}$ and $\Sigma_{3}$ with $\mbox{meas}(\Sigma_{1})>0$. The body is  clamped on  $\Sigma_{1}$  and subjected to the action of volume \mbox{force with  density} $\bm{f}_{0}$.
An unknown surface traction ( for convenience, we denote by $\bm{f}_{2}$  its density) with impulsive effect  is applied on $\Sigma_{2}$. On $\Sigma_{3}$, the body may contact with an obstacle.  We do not show expressly the relation of various functions and $\bm{y}$.

Let $\bm{\nu}$ be unit outward normal  vector,  $\mathbb{S}^{n}$ be the space of symmetric matrix of order two on $R^{n}$. $\mathbb{S}^{n}$ and $R^{n}$ are equipped with, respectively, the following inner products and  norms:
\begin{align*}
& \bm{\xi\cdot \zeta}=\xi_{ij}\zeta_{ij},\qquad \|\bm{\xi}\| =(\bm{\xi\cdot \xi})^{\frac{1}{2}}, \qquad\forall   \bm{\xi}  ,\bm{\zeta}  \in \mathbb{S}^{n}.\\
&\bm{m\cdot n}=m_{i}n_{i},\qquad\;\; \|\bm{m}\| =(\bm{m\cdot m})^{\frac{1}{2}}, \qquad\forall   \bm{m}  ,\bm{n}  \in R^{n}.
\end{align*}
Here, the summation convention  is adopted. For any $\bm{\eta}\in R^{n}$ and $\bm{\sigma}\in \mathbb{S}^{n}$, we denote by $\eta_{\nu}=\bm{\eta}\cdot\bm{\nu}$ the normal components of $\bm{\eta}$, $\bm{\eta}_{\tau}=\bm{\eta}-\eta_{\nu}\bm{\nu}$ the tangential components of $\bm{\eta}$, $\sigma_{\nu}$=$(\bm{\sigma\nu})$$\cdot\bm{\nu}$ the normal components of $\bm{\sigma}$,  $\sigma_{\nu}$=$(\bm{\sigma\nu})$$\cdot\bm{\nu}$ the  components of $\bm{\sigma}$.

We also denote by $\bm{u}=(u_i)\in R^{n}$, $\bm{\sigma}\in \mathbb{S}^{n}$ and $\bm{\varepsilon(u)}=(\varepsilon_{ij}(\bm{u}))\in \mathbb{S}^{n}$, respectively,  the displacement vector,  the stress tensor and the linearized (small) strain tensor, where
$$ \varepsilon_{ij}(\bm{u})=\frac{1}{2}\left(u_{i,j}+u_{j,i} \right),\quad u_{i,j}=\frac{\partial u_{i}}{\partial y_{i}},   \quad \bm{y}=(y_{i})\in V\cup \partial V,  \quad i,j=1,\cdots,n.$$
For more details, we refer the reader to \cite{Han2017,bh2013,SF2012, narwaa, Zeng3}.

 We now turn to present a new contact problem with the surface traction governed by a fractional impulsive differential equation.

\begin{problem}\label{p2} Find a stress $\bm{\sigma}:V\times Q\to \mathbb{S}^{n}$, a surface traction density $\bm{f}_{2}:\Sigma_2\times Q\to R^{n}$ and a displacement field $\bm{u}:V\times Q\to  R^{n}$  such that
\begin{align}
\bm{\sigma}(t)=\mathbb{A}\bm{\varepsilon}(\bm{u}(t))\quad&\quad \mbox{in}\quad \;V\times Q,\label{e6b.1}\\
\mbox{Div}\bm{\sigma}(t)+\bm{f}_{0}(t)=\bm{0}\quad&\quad \mbox{in}\quad \;V\times Q,\label{e6n.2}\\
\bm{u}(t)=\bm{0}\quad&\quad \mbox{on}\quad \Sigma_{1}\times Q,\label{e6.3}\\
\bm{\sigma}(t)\bm{\nu}=\bm{f}_{2}(t)\quad&\quad \mbox{on}\quad \Sigma_{2}\times Q,\label{e6.4}\\
_{0} ^{C}D^{\kappa} _{t}\bm{f}_{2}(t)=F(t,\bm{f}_{2}(t),\bm{u}(t)) \quad&\quad \mbox{on}\quad \;\Sigma_{2}\times Q,\label{e6.8}\\
t\in Q,\;0<\kappa<1,\;t\neq \tau_j,\;j=1,2,\cdots,m\quad&\label{e6e.z8}\\
\Lambda \bm{f}_{2}(\tau_j)=\Theta_j(\bm{f}_{2}(\tau_j ^{-})),\;j=1,2,\cdots,m \quad&\quad \mbox{on}\quad \;\Sigma_{2}\times Q,\label{e6.08}\\
\bm{f}_{2}(0)=\bm{f}_{2} ^{0}\quad&\quad \mbox{on}\quad \;\Sigma_{2}\times Q,\label{e6.09}\\
-\bm{\sigma}_{\tau}(t)\in\partial j_{\tau}(\bm{u}_{\tau}(t))\quad&\quad \mbox{on}\quad \;\Sigma_{3}\times Q,\label{e6.7}\\
-\sigma_{\nu}(t)\in \partial j_{\nu}(u_{\nu}(t))\quad&\quad\mbox{on}\quad \Sigma_{3}\times Q. \label{e6.5}
\end{align}
Here, relation (\ref{e6b.1}) presents an elastic \mbox{constitutive law} with $\mathbb{A}$ being the elasticity operator. Equation (\ref{e6n.2})  \mbox{is  the equilibrium equation} and equation \eqref{e6.3} implies that the body is \mbox{clamped on} $\Sigma_{1}$. The equalities (\ref{e6.4})-(\ref{e6.09}) show that the traction is acted on $\Sigma_{2}$ and the density of the surface traction is \mbox{governed by} a fractional impulsive differential equation, where the $F$ is a function  to be specified later. The set-valued relations \eqref{e6.7} and \eqref{e6.5} denote,  the friction and contact conditions, respectively, where $j_{\tau}$ and $j_{\nu}$ are  locally Lipschitz functionals.
\end{problem}

To deduce  the weak   formulation of Problem \ref{p2}, we consider spaces
$$\mathcal{H}=L^{2}(V;\mathbb{S}^{n})^{n\times n},\quad \mathcal{V}=\{\bm{v}\in H^{1}(V;R^{n})| \bm{v}=\bm{0}\;\mbox{on}\;\Sigma_{1}\}$$
equipped with  \mbox{the inner} products
$$(\bm{\sigma},\bm{\tau})_{\mathcal{H}}=\int_{V}\sigma_{ij}\tau_{ij}dx,\quad (\bm{u},\bm{v})_{\mathcal{V}}=(\bm{\varepsilon}(\bm{u}),\bm{\varepsilon}(\bm{v})_{\mathcal{H}}$$
and corresponding norms \mbox{}$\|\cdot\|_{\mathcal{H}}$ and $\|\cdot\|_{\mathcal{V}}$, respectively.
We denote by $\mathcal{V}^{*}$ the dual space of $\mathcal{V}$, $\langle\cdot,\cdot\rangle_{\mathcal{V}^{*}\times \mathcal{V}}$ the  duality pairing between $\mathcal{V}^{*}$ and $\mathcal{V}$. The trace theorem states
\begin{align}\label{e.6.v11}
\|\gamma{\bm{v}}\| _{L^{2}(\Sigma_{3};R^{n})}\leq \|\gamma\|\|\bm{v}\|_{\mathcal{V}},\quad \forall \bm{v}\in \mathcal{V},
\end{align}
where $\gamma$ is the trace operator defined by
$$\gamma:\mathcal{V}\to L^{2}(\Sigma_3;R^{n}).$$

 Finally, we formulate the Green's formula,  which will be used in the rest of the paper.
\begin{align}\label{ee531}
\int_{V}\bm{\sigma\cdot\varepsilon(u)}dy+\int_{V}\mbox{Div}\bm{\sigma\cdot u}dy=\int_{\partial V}\bm{\sigma\nu\cdot u}d\tau.
\end{align}
In order to study Problem \ref{p2},  we impose some hypotheses on the relevant data.
\begin{itemize}
\item[$\bm{H}(\mathbb{A})$:] The elasticity operator $\mathbb{A}=(\mathbb{A}_{ijkl}):V\times \mathbb{S}^{n}\to \mathbb{S}^{n}$ satisfies the following conditions:
\begin{itemize}
\item[(i)] $\mathbb{A}_{ijkl}=\mathbb{A}_{klij}= \mathbb{A}_{jikl}\in L^{\infty}(V)$;
\item[(ii)] $\mathbb{A}(\bm{y},\bm{0})\in \mathcal{H}$ for a.e. $\bm{y}\in V$;
\item[(iii)] there exists $L_{\mathbb{A}}>0$ such that
$$\|\mathbb{A}(\bm{y,\zeta}_{1})-\mathbb{A}(\bm{y,\zeta}_{2})\|\leq L_{\mathbb{A}}\|\bm{\zeta}_{1}-\bm{\zeta}_{2}\|,\quad \mbox{for all}\; \bm{\zeta}_{1},\bm{\zeta}_{2}\in \mathbb{S}^{n}, \mbox{a.e.}\; \bm{y}\in V;$$
\item[(iv)]there exists $m_{\mathbb{A}}>0$ such that
$$(\mathbb{A}(\bm{y,\zeta}_{1})-\mathbb{A}(\bm{y,\zeta}_{2}))\cdot(\bm{\zeta}_{1}-\bm{\zeta}_{2})\geq m_{\mathbb{A}}\|\bm{\zeta}_{1}-\bm{\zeta}_{2}\|^{2},\quad \forall \bm{\zeta}_{1},\bm{\zeta}_{2}\in \mathbb{S}^{n}.$$
\end{itemize}
\item[$\bm{H}(F)$:] The  function  $F:Q\times \Sigma_2\times L^{2}(\Sigma_2;R^{n})\times \mathcal{V}\to L^{2}(\Sigma_2;R^{n})$  is such that
\begin{itemize}
\item[(i)]$F(\cdot,\bm{x},\bm{y},\bm{z})$ is continuous for all $(\bm{y},\bm{z})\in L^{2}(\Sigma_2;R^{n})\times \mathcal{V}$, a.e. $\bm{x}\in \Sigma_2$;
\item[(ii)]there exists $M_1>0$ such that
$$\|F(t,\bm{x},\bm{z}_{1},\bm{y}_{1})-F(t,\bm{x},\bm{z}_{2},\bm{y}_{2})\|\leq M_1\left(\|\bm{z}_{1}-\bm{z}_{2}\|+\|\bm{y}_{1}-\bm{y}_{2}\|\right), \; $$
for all $(t,\bm{z}_{i},\bm{y}_{i})\in Q\times L^{2}(\Sigma_2;R^{n})\times \mathcal{V} (i=1,2)$, a.e. $\bm{x}\in \Sigma_2$;
\item[(iii)] there exists $\phi\in  L_{+} ^{\frac{1}{p}}[0,T](0<p<\kappa<1)$ satisfying
$$\|F(t,\bm{x},\bm{z},\bm{y})\|\leq \phi(t),\;\mbox{for all}\; (t,\bm{z},\bm{y})\in Q \times L^{2}(\Sigma_2;R^{n})\times \mathcal{V},\; \mbox{a.e.}\; \bm{x}\in \Sigma_2.$$
\end{itemize}

\item[$\bm{H}(I)$:]  $\Theta_j:L^{2}(\Sigma_2;R^{n})\to L^{2}(\Sigma_2;R^{n})(j= 1,2,...,m)$ is \mbox{bounded} and there \mbox{exist}  $d_j> 0$ satisfying
$$\|\Theta_j(\bm{z}_1)-\Theta_j(\bm{z}_1)\|_{L^{2}(\Sigma_2;R^{n})}\leq d_j\|\bm{z}_1-\bm{z}_2\|_{L^{2}(\Sigma_2;R^{n})},\quad \forall \bm{z}_1,\bm{z}_2\in L^{2}(\Sigma_2;R^{n}).$$

\item[$\bm{H}(j_{\nu})$:] The  function $j_{\nu}:\Sigma_{3}\times R\to R$  is such that
\begin{itemize}
\item[(i)] For  a.e. $\bm{y}\in \Sigma_{3}$, $j_{\nu}(\bm{y},\cdot)$ is locally Lipschitz on $R$;
\item[(ii)] For all $r\in R$, $j_{\nu}(\cdot,r)$ is measurable on $\Sigma_{3}$ ;
\item[(iii)] For all $r\in R$, a.e. $\bm{y}\in \Sigma_{3}$, there exist $\overline{c}_{0}\geq 0$ such that
$$|\partial j_{\nu}(\bm{y},r)|\leq \overline{c_{0}}(1+|r|);$$
\item[(iv)] For all $s_{i}\in R(i=1,2)$, a.e. $\bm{y}\in \Sigma_{3}$, there exist $\alpha_{\nu1}>0$ such that
\begin{align*}
 j_{\nu} ^{\circ}(\bm{y},s_{1};s_{2}-s_{1})+j_{\nu} ^{\circ}(\bm{y},s_{2};s_{1}-s_{2})\leq \alpha_{\nu1}|s_{1}-s_{2}|^{2}.
 \end{align*}
\end{itemize}

\item[$\bm{H}(j_{\tau})$:] The  function $j_{\tau}:\Sigma_{3}\times R^{n}\to R$  is such that
\begin{itemize}
\item[(i)]  $j_{\tau}(\bm{y},\cdot)$ is locally Lipschitz on $R$ for all  a.e. $\bm{y}\in \Sigma_{3}$;
\item[(ii)]$j_{\tau}(\cdot,\bm{r})$ is measurable on $\Sigma_{3}$ for all $\bm{r}\in R^{n}$;
\item[(iii)] there exist $\overline{c}_{1}\geq 0$ such that
$$|\partial j_{\tau}(\bm{y},\bm{r})|\leq \overline{c_{1}}(1+\|\bm{r}\|),$$
 for all $\bm{r}\in R^{n}$, a.e. $\bm{y}\in \Sigma_{3}$;
\item[(iv)]there exist $\alpha_{\nu2}>0$ such that
\begin{align*}
 j_{\tau} ^{\circ}(\bm{y},\bm{s}_{1};\bm{s}_{2}-\bm{s}_{1})+j_{\tau} ^{\circ}(\bm{y},\bm{s}_{2};\bm{s}_{1}-\bm{s}_{2})\leq \alpha_{\nu2}\|\bm{s}_{1}-\bm{s}_{2}\|^{2},
\end{align*}
for all $\bm{s}_{i}\in R(i=1,2)$ and a.e. $\bm{y}\in \Sigma_{3}$.
\end{itemize}

\item[$\bm{H}(\bm{f})$:] For the densities of \mbox{}body force satisfies that
$$\bm{f}_{0}\in \mathcal{I}C(Q;L^{2}(V;R^{n})).$$
\item[$\bm{H}(0)$:]
\begin{itemize}
\item[(i)]$m_{\mathbb{A}}>(\alpha_{\nu1}+\alpha_{\nu2})c_0 ^{2}$;
\item[(ii)] $\frac{T^{\kappa}M_1 c_0}{\kappa\left[m_{\mathbb{A}}- (\alpha_{\nu1}+\alpha_{\nu2})c_0 ^{2}\right]\Gamma(\kappa)}<1.$
\end{itemize}
\end{itemize}
Utilizing the Green formula \eqref{ee531}, we get the variational form of Problem \ref{p2}.
\begin{problem}\label{p3} Find a displacement field $\bm{u}:Q\to \mathcal{V}$ and  a  surface traction density $\bm{f}_{2}: Q\to L^{2}(\Sigma_2;R^{n})$ such that
$$
\begin{cases}
_{0} ^{C}D^{\kappa} _{t}\bm{f}_{2}(t)=F(t,\bm{f}_{2}(t),\bm{u}(t)),\quad t\in Q,\;0<\kappa<1,\;t\neq \tau_j,\;j=1,2,\cdots,m,\\
\Lambda \bm{f}_{2}(\tau_j)=\Theta_j(\bm{f}_{2}(\tau_j ^{-})),\;j=1,2,\cdots,m,\\
\bm{f}_{2}(0)=\bm{f}_{2} ^{0},\\
\left(\mathbb{A}\bm{\varepsilon}(\bm{u}(t)),\bm{\varepsilon}(\bm{v})\right)_{\mathcal{H}}
+\int_{L_{3}}\left(j_{\nu} ^{\circ}(u_{\nu}(t);v_{\nu})+j_{\tau} ^{\circ}(\bm{u}_{\tau}(t);\bm{v}_{\tau})\right)da\notag\\
\quad \geq \int_{L_{2}}\bm{f}_{2}(t)\bm{v}da+\int_{V}\bm{f}_{0}(t)\bm{v}d\bm{x},\quad \forall (t,\bm{v})\in Q\times\mathcal{V}.
\end{cases}
$$
\end{problem}

\subsection{ Existence and uniqueness for the contact problem}
We  define the maps $A:\mathcal{V}\to \mathcal{V}^{*}$, $f:Q\times L^{2}(\Sigma_2;R^{n}) \times \mathcal{V}\to L^{2}(\Sigma_2;R^{n})$, $J:L^{2}(\Sigma_{3};R^{n})\to R$, and  $g:L^{2}(\Sigma_2;R^{n})\to \mathcal{V}^{*}$ by setting
\begin{align}
&\langle A\bm{u},\bm{v}\rangle_{\mathcal{V}^{*}\times\mathcal{V}}=\left(\mathcal{E}\bm{\varepsilon (u)},\bm{\varepsilon (v)}\right)_{\mathcal{H}},\label{ee.6.16}\\
&f(t,\bm{f}_2,\bm{v})=F(t,\bm{f}_{2}(t),\bm{v}(t)),\\
&J(\gamma\bm{u})=\int_{\Sigma_{3}}\left(j_{\nu} (u_{\nu}(t))+j_{\tau} (\bm{u}_{\tau}(t))\right)da,\label{ee.6.0.19}\\
&\langle g(\bm{f}_{2}(t)),\bm{v}\rangle_{\mathcal{V}^{*}\times\mathcal{V}}=\int_{V}\bm{f}_{0}(t)\cdot\bm{v}d\bm{x}+\int_{\Sigma_{2}}\bm{f}_{2}(t)\cdot\bm{v}da,\label{ee.6.20}
\end{align}
for all $(t,\bm{f}_{2},\bm{u},\bm{v})\in Q\times R^{n}\times \mathcal{V}\times \mathcal{V}$.

Then,  Problem \ref{p3} is equivalent to the  problem:
\begin{problem}\label{p4}
 Find a displacement vector $\bm{u}:Q\to \mathcal{V}$ and a  surface traction density $\bm{f}_{2}: Q\to L^{2}(\Sigma_2;R^{n})$ such that
$$
\begin{cases}
_{0} ^{C}D^{\kappa} _{t}\bm{f}_{2}(t)=f(t,\bm{f}_{2}(t),\bm{u}(t)),\;t\in Q,\;0<\kappa<1,\;t\neq \tau_j,\;j=1,2,\cdots,m,\\
\Lambda \bm{f}_{2}(\tau_j)=\Theta_j(\bm{f}_{2}(\tau_j ^{-})),\;j=1,2,\cdots,m,\\
\bm{f}_{2}(0)=\bm{f}_{2} ^{0},\\
\langle A\bm{u}(t),\bm{v}\rangle_{\mathcal{V}^{*}\times\mathcal{V}}+J^{\circ}(\gamma\bm{u};\gamma\bm{v})\geq \langle g(\bm{f}_{2}(t)),\bm{v}\rangle_{\mathcal{V}^{*}\times\mathcal{V}},
\quad \forall (t,\bm{v})\in t\times\mathcal{V}.
\end{cases}
$$
\end{problem}
Clearly, Problem \ref{p4} is the form of  Problem \ref{p1} with $Z_{1}= L^{2}(\Sigma_2;R^{n}),Z_2=\mathcal{V},Y=L^{2}(\Sigma_3;R^{n})$.

\begin{theorem}\label{hvinqdl}
Problem \ref{p4} admits a \mbox{unique} solution $(\bm{f}_{2},\bm{u}(t))\in \mathcal{I}C(Q;L^{2}(\Sigma_2;R^{n})) \times \mathcal{I}C(Q;\mathcal{V})$ providing that hypotheses $\bm{H}(\mathbb{A})$, $\bm{H}(F)$,$\bm{H}(I)$, $\bm{H}(j_{\nu})$, $\bm{H}(j_{\tau})$, $\bm{H}(\bm{f})$ and $\bm{H}(0)$  hold.
 \end{theorem}

\textbf{Proof.} To prove  Theorem \ref{hvinqdl}, we only need  to check the validity of assumptions (HA), (Hf), (HI), (HN), (HJ), (Hg), and (HO).

Firstly, conditions $\bm{H}(\mathbb{A})$, $\bm{H}(F)$ and $\bm{H}(I)$ indicate that assumptions (HA), (Hf) and (HI) are fulfilled with $m_{A}=m_{\mathbb{A}}$. Since the the trace operator is compact and surjective, we see that the assumption (HN) holds.

Clearly, \eqref{ee.6.20} implies that the assumption (Hg) holds with $m_g=\|\gamma\|$. By hypotheses $\bm{H}(j_{\nu})$, $\bm{H}(j_{\tau})$ and Lemma 14 in \cite{Rothe2019}, it follows from Lemma 14 in \cite{Rothe2019} that the functional $J$ in \eqref{ee.6.0.19} is locally Lipschitz on $L^{2}(\Sigma_{3};R^{n})$ and
\begin{align}\label{j545.e}
  J^{\circ}(\gamma\bm{u};\gamma\bm{w})=\int_{L_{3}}\left(j^{\circ} _{\nu} (u_{\nu}(t);w_{\nu})+j^{\circ} _{\tau} (\bm{u}_{\tau}(t);\bm{w}_{\tau})\right)da,\quad \forall \bm{u},\bm{w}\in \mathcal{V}
\end{align}
is the \mbox{generalized} directional \mbox{derivative} of $J$ at $\gamma\bm{u}$ in \mbox{the} directional $\gamma\bm{w}$.

Moreover, the assumption (HJ) holds with $c_J=\alpha_{\nu1}+\alpha_{\nu2}$ and $m_J=\max\{\overline{c_0},\overline{c_1}\}$.
Combining Theorem \ref{hv+zxdl} with the hypothesis $\bm{H}(0)$,  we see that Theorem \ref{hvinqdl} holds. \hfill$\Box$

\subsection{ A  convergence result for the contact problem}

The above analysis reveals that the solution  of Problem \ref{p4} relies on the data $j_{\nu}$ and $j_{\tau}$. In what follows,  we present a  continuous \mbox{}dependence result of the solution in relation to these data. \mbox{}We consider the  perturbation data $j_{\nu\delta}$ and $j_{\tau\delta}$ of $j_{\nu}$ and $j_{\tau}$, respectively, which satisfy hypotheses $\bm{H}(j_{\nu})$ and $\bm{H}(j_{\tau})$.  For each $\delta>0$, define a function
 $J_{\delta}:L^{2}(\Sigma_{3};R^{n}) \to R$  by setting
\begin{align}\label{J547as}
  J_{\delta}(\gamma\bm{u})=\int_{\Sigma_{3}}\left(j_{\nu\delta} (u_{\nu}(t))+j_{\tau\delta} (\bm{u}_{\tau}(t))\right)da,\quad \forall \bm{u}\in\mathcal{V}.
\end{align}

The perturbation  problem \mbox{of} Problem \ref{p4} can be formulated as follows.
\begin{problem}\label{p5}
 Find a displacement vector $\bm{u}_{\delta}:Q\to \mathcal{V}$ and a  surface traction density $\bm{f}_{2\delta}: Q\to L^{2}(\Sigma_2;R^{n})$ such that
\begin{align}
&_{0} ^{C}D^{\kappa} _{t}\bm{f}_{2\delta}(t)=f(t,\bm{f}_{2\delta}(t),\bm{u}_{\delta}(t)),\label{e16vc.08}\\
&t\in Q,\;0<\kappa<1,\;t\neq \tau_j,\;j=1,2,\cdots,m,\\
&\Lambda \bm{f}_{2\delta}(\tau_j)=\Theta_j(\bm{f}_{2\delta}(\tau_j ^{-})),\;j=1,2,\cdots,m,\\
&\bm{f}_{2\delta}(0)=\bm{f}_{2} ^{0},\\
&\langle A\bm{u}_{\delta}(t),\bm{v}\rangle+J^{\circ} _{\delta}(\gamma\bm{u}_{\delta};\gamma\bm{v})\geq \langle g(\bm{f}_{2\delta}(t)),\bm{v}\rangle,
\quad \forall (t,\bm{v})\in t\times\mathcal{V}.\label{e6bbe.23}
\end{align}
\end{problem}

 Denote the constants involved in hypotheses $\bm{H}(j_{\nu\delta})$(iv) and $\bm{H}(j_{\tau\delta})$(iv) by $\alpha_{\nu1\delta}$ and $\alpha_{\nu2\delta}$, respectively.  In addition, we impose the following hypotheses on the data.
\begin{itemize}
\item[$\bm{H}(j^{*})$:]There \mbox{}exists a function $\overline{V}:R^{+}\to R^{+}$ satisfying
\begin{itemize}
\item[(i)] $|\partial j_{\nu} (\bm{x},r)-\partial j_{\nu\delta}(\bm{x},r)|\leq \overline{V}(\delta)|r|, \; \mbox{for all} \;(\delta,r)\in R^{+}\times R$, \mbox{a.e.} $\bm{x}\in  \Sigma_{3}$;
\item[(ii)] $\|\partial j_{\tau} (\bm{x},\bm{b})-\partial j_{\tau\delta} (\bm{x},\bm{b})\|\leq \overline{V}(\delta)\|\bm{b}\|, \; \forall \;(\bm{x},\bm{b})\in \Sigma_3\times R^{n};$
\item[(iii)]$\lim_{\delta\to 0} \overline{V}(\delta)=0$.
\end{itemize}

\item[$\bm{H}(0^{*})$:] There \mbox{}exists $m_{\mathbb{A}_{0}}>0$ \mbox{}such that
\begin{itemize}
\item[(i)]$m_{\mathbb{A}}>m_{\mathbb{A}_{0}}>(\alpha_{\nu1\delta}+\alpha_{\nu2\delta})c_{0} ^{2}$;
\item[(ii)] $\frac{T^{\kappa}M_1 c_{0}}{\kappa\left[m_{\mathbb{A}}- (\alpha_{\nu1\delta}+\alpha_{\nu2\delta})c_{0} ^{2}\right]\Gamma(\kappa)}<1.$
\end{itemize}
\end{itemize}

\begin{theorem}\label{hv+q5dl}
Assume that hypotheses $\bm{H}(\mathbb{A})$, $\bm{H}(F)$, $\bm{H}(I)$, $\bm{H}(j_{\nu})$, $\bm{H}(j_{\tau})$, $\bm{H}(\bm{f})$, $\bm{H}(j^{*})$,  $\bm{H}(0)$ and $\bm{H}(0^{*})$ hold. Then
\begin{itemize}
\item[(a)] for each $\delta>0$, Problem \ref{p5} has a unique solution $(\bm{f}_{2\delta},\bm{u}_{\delta}(t))\in \mathcal{I}C(Q;L^{2}(\Sigma_2;R^{n})) \times \mathcal{I}C(Q;\mathcal{V})$.
\item[(b)]$(\bm{f}_{2\delta},\bm{u}_{\delta}(t))$ converges to $(\bm{f}_{2},\bm{u}(t))$, the solution  of Problem \ref{p4}, i.e.,
\begin{align}\label{e4v.cd070a}
(\bm{f}_{2\delta},\bm{u}_{\delta}(t))\to (\bm{f}_{2},\bm{u}(t))\quad \mbox{as}\;\delta\to 0 \;\mbox{for all}\;t\in Q.
\end{align}
\end{itemize}
 \end{theorem}
\textbf{Proof.} (a) In view of Theorem \ref{hvinqdl}, the proof is obvious.

(b) We employ Theorem \ref{hvpre5dl} to prove \eqref{e4v.cd070a}. To this end, we only need to check the validity of  assumptions (HO$^{*}$) and (HJ$^{*}$).
Clearly, the hypothesis (H0$^{*}$) implies that the assumption  (HO$^{*}$) holds.  By Proposition 3.35 of \cite{bh2013}, Corollary 4.15 \cite{bh2013} and the hypothesis $\bm{H}(j^{*})$, for any $(\bm{u},\bm{\xi},\bm{\xi}_{\delta})\in \mathcal{V}\times\gamma^{*} \partial J(\gamma\bm{u})\times \gamma^{*} \partial J_{\delta}(\gamma\bm{u})$ and $(\xi_{\nu},\xi_{\nu\delta},\bm{\xi}_{\tau},\bm{\xi}_{\tau\delta})\in \partial j_{\nu}(u_{\nu}(t))\times\partial j_{\nu\delta} (u_{\nu}(t))\times \partial j_{\tau}(\bm{u}_{\tau}(t))\times \partial j_{\tau\delta}(\bm{u}_{\tau}(t))$,  we have
  \begin{align*}
    \|\bm{\xi}-\bm{\xi}_{\delta}\|&\leq \|\gamma^{*}\|\int_{L_{3}}\left(|\xi_{\nu}-\xi_{\nu\delta}|+\|\bm{\xi}_{\tau}-\bm{\xi}_{\tau\delta}\|\right)da \\
    &\leq  \|\gamma^{*}\|\overline{V}(\delta)\int_{L_{3}}\left(|u_{\nu}|+\|\bm{u}_{\tau}\|\right)da\\
    &\leq  \left(\|\gamma\|^{2}\mbox{meas}(\Sigma_{3})\|\bm{u}\|\right)\overline{V}(\delta),
 \end{align*}
which shows that the assumption (HJ$^{*}$) holds with $V(\delta)=\left(\|\gamma\|^{2}\mbox{meas}(\Sigma_{3})\|\bm{u}\|\right)\overline{V}(\delta)$. The convergence result (\ref{e4v.cd070a})  now follows from Theorem \ref{hvpre5dl}.  \hfill  $\Box$

\end{document}